\documentclass[11pt, reqno]{amsart}
\usepackage{amscd, amsmath, amsthm, amssymb}

\usepackage[bookmarksnumbered,colorlinks, plainpages]{hyperref}
\usepackage{stackrel}
\hypersetup{
pdftitle={Research Paper},
pdfauthor={Reza Abdolmaleki},
pdfsubject={on the Betti number of monomial ideals and their powers},
pdfcreator={Reza Abdolmaleki},
colorlinks,
linkcolor=blue,
citecolor=magenta,
anchorcolor=red,
bookmarksopen,
urlcolor=red,
filecolor=red,
}
\theoremstyle{plain}
\newtheorem{Theorem}{Theorem}[section]

\newtheorem{Corollary}[Theorem]{Corollary}
\newtheorem{Proposition}[Theorem]{Proposition}

\theoremstyle{definition}
\newtheorem{Example}{Example}

\newtheorem{Notation}{Notation}

\theoremstyle{remark}
\newtheorem{Remark}{Remark}

\newcommand{\KK}{\mathbb{K}}
\newcommand{\ZZ}{\mathbb{Z}}

\newcommand{\height}{\mathrm{height}}

\newcommand{\ab}{\mathbf{a}}
\newcommand{\cb}{\mathbf{c}}
\newcommand{\eb}{\mathbf{e}}
\newcommand{\bb}{\mathbf{b}}

\newcommand{\xb}{\mathbf{x}}
\def\cocoa{{\hbox{\rm C\kern-.13em o\kern-.07em C\kern-.13em o\kern-.15em A}}}

\begin{document}
\author[R. Abdolmaleki]{Reza Abdolmaleki}
\email{reza.abd110@gmail.com, abdolmaleki@iasbs.ac.ir}
\address{Reza Abdolmaleki, School of Mathematics, Institute for Research in Fundamental Sciences (IPM), P.O. Box: 19395-5746, Tehran, Iran}
\address{Department of Mathematics, Institute for Advanced Studies in Basic Sciences (IASBS), Zanjan 45137-66731, Iran}

\author[R.  Zaare-Nahandi]{{Rashid} {Zaare-Nahandi}}
\email{rashidzn@iasbs.ac.ir}
\address{Rashid Zaare-Nahandi, Department of Mathematics, Institute for Advanced Studies in Basic Sciences (IASBS), Zanjan 45137-66731, Iran}

\title{On the Betti numbers of monomial ideals and their powers}

\begin{abstract}
Let $S=\KK[x_1,\ldots,x_n]$ the polynomial ring  over a field $\KK$.  In this paper for some families of monomial ideals $I \subset S$  we  study the minimal number of generators of $I^k$. We use  this results  to find some other Betti numbers of these families of ideals for special choices of $n$, the number of variables.
\end{abstract}

\thanks{This research was in part supported by a grant from IPM ((No.1401130016). The first author would likes to thank the Institute for Research in Fundamental Sciences (IPM), Tehran, Iran.}

\subjclass[2020]{Primary: 13F20; Secondary: 05E40.}
\keywords{Monomial ideal, Betti number, $t$-spread monomial,  ideal of Veronese type}

\maketitle

\section{Introduction}

Using the structure of an ideal $I$ in a commutative ring to find the Betti numbers of  $I$  and the powers $I^k$ is a complicated problem. In paricular, finding $\mu(I)$, the minimal number of generators of a graded  polynomial ideal  $I$ and predicting the beahaviour of the function $\mu(I^k)$ is quite difficult and has been studied a lot (for instance, see \cite{AHZN}, \cite{AK}, \cite{E1}, \cite{HHZ}, \cite{HQS}, \cite{HSZ} and \cite{HZ1}). In this paper we find the minimal generators of some families of equigenerated monomial ideals (monomial ideals generated in a single degree)  in the polynomial ring $\KK[x_1, \ldots, x_n]$ over a field $\KK$.  Moreover, we find some other Betti numbers of these ideals for special choices of $n$, the number of variables.

Let $\KK$ be a field and $S=\KK[x_1,\ldots,x_n]$ be the polynomial ring in the variables $x_1,\ldots,x_n$ over $\KK$ . Also, let $I$ be a graded ideal in $S$ and 
\[
0\to S^{\beta_n}\to \cdots \to S^{\beta_2}\to S^{\beta_1}\to S \to S/I\to  0
\]
be the minimal free resolution of $S/I$. The numbers $\beta_1, \ldots, \beta_n$ are called the Betti numbers of $S/I$.

An equigenerated monomial ideal $I$ with the minimal set of generators $G(I)$ is called a {\em polymatroidal ideal} if for any pairs of monomial $x_1^{a_1} \ldots x_n^{a_n}$ and $x_1^{a'_1} \ldots x_n^{a'_n}$  in $G(I)$ with the property that $a_i>a'_i$ for some $i$, there exists a $j$ such that $a_j<a'_j$ and $(x_j/x_i)(x_1^{a_1} \ldots x_n^{a_n}) \in G(I)$.
We say that the ideal $I$ has a $d$-linear resolution if the graded minimal free resolution of $S/I^k$ is of the form 
\[
0\to S(-d-s)^{\beta_s}\to \cdots \to S(-d-1)^{\beta_2}\to S(-d)^{\beta_1}\to S\to S/I^k\to  0.
\]

Let $I$ be a polymatroidal ideal. Since all powers of a  polymatroidal ideal are polymatroidal (\cite[Theorem12.6.3]{HH})  and  a polymatroidal ideal have a  linear resolution (\cite{HT}), the minimal free resolution of  $S/I^k$  is of the form
  \[
0\to S(-(kd +n-1))^{\beta_n^k}\to \cdots \to S(-(kd+1))^{\beta_2^k}\to S(-kd)^{\beta_1^{k}}\to S\to S/I^k\to  0.
\]
where $\beta_i^{k}=\beta_i(S/I^k)$.

An important class of polymatroidal ideals is the class of ideals of Veronese type. Fix  integer $d$ and the integer vector  $\ab=(a_1, \ldots , a_n)$  with $d \geq a_1 \geq \ldots \geq a_n \geq 1$. An {\em ideal of Veronese type} is an  ideal $I_{\ab,n,d}$ with the following minimal set of generators
\[
G(I_{\ab,n,d})=\{x_1^{b_1}x_2^{b_2}\cdots x_n^{b_n}\; \mid \; \sum_{i=1}^nb_i=d \text{ and  $b_i\leq a_i$ for $i=1,\ldots,n$}\}.
\]

In Section \ref{1} we find the minimal set of generators of all powers of ideals of Veronese type. Also, we use this result to find the minimal number of generators ($\beta_1$) of some other classes of equigenerated monomial ideals. In Section \ref{2} we use $\beta_1$ to find some other Betti numbers of these families of  ideals for special choices of $n$.

\section{The minimal number of genrators of some monomial ideals}
\label{1}
 Let  $n, d \geq 1$  and $t \geq 0$ be fixed  integers. The following notations are obtained from \cite{DHQ}.  We denote by $\mathcal{A}_{n,d}$ the set of all multisets $ A \subset [n]$ with $|A|=d$. A multiset $ \{i_1 \leq  i_2 \leq  \ldots \leq  i_d\} \subset [n]$ is called {\em $t$-spread}, if $ i_{j+1}- i_{j}\geq t$ for all $j$. The set of all $t$-spread multisets in $\mathcal{A}_{n,d}$ is denoted by $\mathcal{A}_{n,d,t}$. Let $A \subset \mathcal{A}_{n,d,t}$ be a $t$-spread multiset. A subset $B \subset A$ is called a {\em block of size $q$}, if  $B= \{i_j, i_{j+1}, \ldots , i_{j+q-1}\} $  with $ i_{k+1}- i_{k}=t$ for all $k$. Let $c$ be a positive integer. The set of all multisets $A \subset \mathcal{A}_{n,d,t}$  such that $|B| \leq c$ for each block $B \subset A$,  is denoted by $\mathcal{A}_{c,(n,d,t)}$.

Let $\ab=(a_1, \ldots , a_n)$ be a vector of integers such  $d \geq a_1 \geq \ldots \geq a_n \geq 1$. For the integer vector $\cb=(c_1, \ldots , c_n)$  we write $\cb \leq \ab$ if $ c_i \leq a_i$ for all $i$.

Let $S=\KK[x_1,\ldots,x_n]$ be the polynomial ring in the variables $x_1,\ldots,x_n$ over a field  $\KK$.  We fix some notations for the following classes of monomial ideals: 

\begin{Notation}
\label{not}

\item[$\bullet$]  We denote by $I_{\ab,n,d}$ the ideal generated by all monomials of degree $d$ whose exponent vectors are boundead by $\ab$. In other words,

\[
G(I_{\ab,n,d})=\{x_1^{b_1}x_2^{b_2}\cdots x_n^{b_n}\; \mid \; \sum_{i=1}^nb_i=d \text{ and  $b_i\leq a_i$ for $i=1,\ldots,n$}\}.
\]
$I_{\ab,n,d}$  is called an ideal of Veronese type.

\item[$\bullet$]  $I_{c,(n,d,t)}:=(\xb_{A} \; \mid \; A \in \mathcal{A}_{c,(n,d,t)})$. The ideal $ I_{c,(n,d,t)}$ is called a {\em $c$-bounded $t$-spread Veronese ideal}. Note that $ I_{c,(n,d,0)}=I_{\cb,n,d}$  where $\cb=(c, \ldots , c) \in \ZZ^n$.
\item[$\bullet$]  We denote by $I_{n,d,t}$  the ideal generated by all $t$-spread monomials  in $S$ of degree $d$. The ideal $I_{n,d,t}$  is called a {\em $t$-spread Veronese ideal} of degree $d$ . One can easily see that $I_{n,d,t} =I_{d,(n,d,t)}$. 
\item[$\bullet$]  The ideal generated by all square free monomials of degree $d$ is called a {\em square free Veronese ideal} of degree $d$  and is denoted by $I_{n,d}$. Recall that a  monomial $x_1^{b_1}x_2^{b_2}\cdots x_n^{b_n} \in S$ is called square free, if $b_i\leq 1$ for all $i$. Therefore, $I_{n,d}=I_{\mathbf{e},n,d}$  where $\mathbf{e}=(1, \ldots , 1) \in \ZZ^n$, and hence $I_{n,d}=I_{1,(n,d,0)}$.
\end{Notation}

In this section we use the structure of the ideals introduced in Notation~\ref{not} to compute their minimal number of generators (and their powers). We denote by  $\mu(I)$ the minimal number of generators of a graded ideal $I \in S$.

Let $\ab=(a_1, \ldots , a_n)$ be a vector of integers such that $d \geq a_1 \geq \ldots \geq a_n \geq 1$.  Set $\alpha_{i,0}^k=0$ and $\alpha_{i,l}^k=\sum_{i =1}^l(ka_i+1)$ for $1\leq i \leq n$,  $1\leq l \leq n$ and $k \geq 1$.
\begin{Theorem}
\label{vt}
Let $I=I_{\ab,n,d}$ be an ideal of Veronese type  with $\ab=(a_1, \ldots , a_n)$. Then 
\[
\mu(I^k)=\sum_{j=0}^n\bigg[(-1)^j \sum_{i=1}^{{n \choose j}}{kd+n-1- \alpha_{i,j}^k \choose n-1}\bigg].
\]
for all $ k \geq 1$.
\end{Theorem}
\begin{proof} 
First we prove the assertion for $k=1$.  In the case that $a_1=a_2= \ldots = a_n=d$, the ideal $I$ is the Veronese ideal of $S$ in degree $d$, that is, the ideal generated by all monomials in $S$ of degree $d$. Therefore, $\mu(I)={d+n-1 \choose n-1}$. Now we assume that $a_i <d$ for some $i$. A typical generator of $I$ is in the form $x_1^{b_1}x_2^{b_2}\cdots x_n^{b_n}$ such that $b_1+b_2+ \ldots + b_n=d$ and $b_i \leq a_i$ for all $i$. We must subtract the bad cases $b_i>a_i$. So we subtract the number of solutions of the equation
 $$b_1+b_2+ \ldots + b_{i-1}+(b_i-a_i-1)+ b_{i+1}+ \ldots +b_n=d-a_i-1,$$ 
which equals to ${d+n-1-(a_i+1)\choose n-1}$. Using the inclusion-exclusion principle we get
\begin{eqnarray*}
\mu(I)&=&{d+n-1 \choose n-1}+\sum_{J \subseteq \{1, \ldots , n\}}(-1)^{|J|}{d+n-1-\sum_{i \in J}(a_i+1)\choose n-1}\\&=&
\sum_{j=0}^n\bigg[(-1)^j \sum_{i=1}^{{n \choose j}}{d+n-1-\alpha_{i,l}^1\choose n-1}\bigg].
\end{eqnarray*}

The assetion for $k\geq 2$ follows from the fact  $(I_{\ab,n,d})^k= I_{k\ab,n,kd}$ by \cite[Lemma~5.1]{HRV}.
\end{proof}

\begin{Remark} 
\label{rem}
In the case that $a_1=a_2= \ldots = a_n=c$ for some positive integer $c$ it is easy to check  that
\[
\mu(I^k)=\sum_{j=0}^{\lfloor \frac{kd}{kc+1}\rfloor}(-1)^j {{n \choose j}}{kd+n-1- j(kc+1) \choose n-1}.
\]
\end{Remark}

\begin{Proposition}
\label{cbtv}
Let $I=I_{c,(n,d,t)}$ be a $c$-bounded $t$-spread Veronese ideal. Then 
\[
\mu(I)=\sum_{j=0}^{\lfloor \frac{d}{c+1}\rfloor}(-1)^{j}{n-(d-1)t \choose j}{n-(d-1)(t-1) -j(c+1)\choose d}.
\]
\end{Proposition}
\begin{proof} 
The ideals $I_{c,(n,d,t)}$  and $I_{c,(n-(d-1)t,d,0)}$ have the same Betti numbers by \cite[Corollary~3.5]{DHQ}. On the other hand,   $I_{c,(n-(d-1)t,d,0)}=I_{\cb,(n-(d-1)t,d}$ where $\cb=(c,\ldots , c) \in \ZZ^n$. So, the desired conclusion follows from Remark~\ref{rem}.
\end{proof} 
\begin{Corollary}
\label{tver}
Let $I=I_{n,d,t}$ be a $t$-spread Veronese ideal of degree $d$. Then 
\[
\mu(I)={n-(d-1)(t-1)\choose d}.
\]
\end{Corollary}
\begin{proof} 
Since $I_{n,d,t}=I_{d,(n,d,t)}$,  the assertion results from Proposition~\ref{cbtv}.
\end{proof}
\begin{Remark} 
An alternative proof for Corollary~\ref{tver} is given in \cite[Theorem~2.3 (d)]{EHQ}.
\end{Remark}
\begin{Proposition}
\label{fff}
Let $I=I_{n,d}$ be a square free Veronese ideal of degree $d$. Then 
\[
\mu(I^k)=\sum_{j=0}^{\lfloor \frac{kd}{k+1}\rfloor}(-1)^{j}{n \choose j}{kd+n-1-j(k+1)\choose kd}.
\]
for all $ k \geq 1$.
\end{Proposition}
\begin{proof} 
The desired conclusion results from Remark~\ref{rem}, since $I_{n,d}=I_{\mathbf{e},n,d}$  where $\mathbf{e}=(1, \ldots , 1) \in \ZZ^n$.
\end{proof}
\section{On the other Betti numbers of our ideals and their powers}
\label{2}
In the previous section we computed the minimal number of generators ($\beta_1$) of  ideals of Veronese type  and their powers.  It is well known that, for a monomial ideal $I$  in $k[x_1,x_2]$  generated by $\mu(I)$ elements, one has $\beta_2=\beta_1-1=\mu(I)-1$ (see \cite[Proposition 3.1]{MS}). In this section, using $\beta_1$ we find the other Betti numbers of  ideals of Veronese type and their powers in $K[x_1,x_2,x_3]$. Moreover, for the other classes of monomial ideals which we studied their first Betti number in Section \ref{1}, we find some of their other  Betti numbers for particular choices of $n$.

For a monomial ideal $I \subset S$ we denote by $\dim(I)$, the Krull dimention of $S/I$.  Let $I=I_{\ab,3,d} \subset \KK[x_1,x_2,x_3]$  be an ideal of Veronese type with $\dim(I)=2$. So, $\height(I)=1$.  Since $a_1 \geq a_2 \geq a_3$, there exists a positive integer $d'$ and a Veronese type ideal $J$ with $\dim(J)=1$ such that $I = x_1^{d'}J$. Indeed,
 \[
d'=\max\{\ell:\; x_1^{\ell}|u\ \ \text{for all }u\in G(I)\}
\] 
and  $J=I_{\bb,3,d-d'}$ where $\bb=(a_1-d',a_2,a_3)$. Set $\delta=d-d'$.
\begin{Proposition}
\label{jjj}
Let $I=I_{\ab,3,d} \subset \KK[x_1,x_2,x_3]$ be an ideal of Veronese type with $\ab=(a_1, a_2 , a_3)$.
Then, for $k\geq1$, if $\dim(I)=0$, 
\[
\beta_2(I^k)=(kd)(kd+2),  \quad  \beta_3(I^k)={kd+1 \choose 2}.
\]
If $\dim(I)=1$, 
\[
\beta_2(I^k)=2\beta_1(I^k)-kd-2,  \quad  \beta_3(I^k)=\beta_1(I^k)-kd-1.
\]
If $\dim(I)=2$,  
\[
\beta_2(I^k)=\beta_1(I^k)-k\delta-2,  \quad  \beta_3(I^k)=\beta_1(I^k)-k\delta-1.
\]
\end{Proposition}
\begin{proof}
Sine all powers of $I$ are polymatroidal, the minimal free resolution of $I$ is of the form
\[
0\to S(-dk-2)^{\beta_2^{k}}\to S(-dk-1)^{\beta_2^{k}}\to S(-dk)^{\beta_1^{k}}\to S\to S/I^k\to  0,
\]
where $\beta_i^{k}=\beta_i(S/I^k)$. Therefor,  
if $\dim(I)=0$, then $I$ is Cohen-Macaulay. Usiung \cite[Theorem~4.1.15]{BH} we get
$\beta_2(I^k)=kd(kd+2)$  and $ \beta_3(I^k)=(kd)(kd+1)/2={kd+1 \choose 2}$.

If  $\dim(I)=1$, using \cite[Theorem~3]{N} we get
\[
\begin{pmatrix}
\beta_2^{k}\\ 
\beta_3^{k} \\ 
\end{pmatrix}
=
\begin{pmatrix}
1& 1 \\ 
0 & 1  \\ 
\end{pmatrix}
\begin{pmatrix}
\beta_1^{k}-{kd\choose 0}\\ 
\beta_1^{k}-{kd+1 \choose 1} \\ 
\end{pmatrix}
=
\begin{pmatrix}
2\beta_1^{k}-kd-2\\ 
\beta_1^{k}-kd-1 \\ 
\end{pmatrix}.
\]
If $\dim(I)=2$, we assume that  $I = x_1^{d'}J$ with $\dim(J)=1$ and set $\delta=d-d'$. 
Since $I$ and $J$ have the same Betti numbers, the assertion follows from the previous case.
\end{proof}
\begin{Example}
\begin{enumerate}
\item[(a)] Let $\ab=(2,2,2)$ and $d=2$. Then $I_{\ab,3,2}=(x_1,x_2,x_3)^2 \subset \KK[x_1,x_2,x_3]$ and so, $\dim(I_{\ab,3,2})=0$. Using Theorem~\ref{vt} and  Proposition~\ref{jjj} we get $\beta_1(I_{\ab,3,2})=6$, $\beta_2(I_{\ab,3,2})=8$ and $\beta_3(I_{\ab,3,2})=3$.
\item[(b)] Let $\eb=(1,1,1)$ and $d=1$. Then $I_{\eb,3,2}=I_{3,2}=(x_1x_2, x_1x_3, x_2x_3) \subset \KK[x_1,x_2,x_3]$ and so, $\dim(I_{\eb,3,2})=1$. Using Corollary~\ref{fff} and  Proposition~\ref{jjj} we get $\beta_1(I_{\eb,3,2})=3$, $\beta_2(I_{\eb,3,2})=2$ and $\beta_3(I_{\eb,3,2})=0$.
\item[(c)] Let $\cb=(8,2,1)$ and $d=8$. Then $I_{\cb,3,8}=(x_1^8,x_1^7x_2, x_1^7x_3,x_1^6 x_2x_3,x_1^6 x_2^2,x_1^5 x_2^2x_3) \subset \KK[x_1,x_2,x_3]$ and so, $\dim(I_{\cb,3,8})=2$. Note that $I_{\cb,3,8}=x_1^5I_{\bb,3,3}$ where $\bb=(3,2,1)$ and hence $\delta=8-5=3$.  Using Theorem~\ref{vt} and  Proposition~\ref{jjj} we get $\beta_1(I_{\cb,3,3})=6$, $\beta_2(I_{\cb,3,3})=7$ and $\beta_3(I_{\cb,3,3})=2$.
\end{enumerate}
\end{Example}
Let $I_{c,(n,d,t)}$ be a $c$-bounded $t$-spread Veronese ideal with $\dim (I)=2$ such that $n-(d-1)t=3$. By \cite[Corollary~3.5]{DHQ} we have $\beta_i(I_{c,(n,d,t)})=\beta_i(I_{c,n-(d-1)t,d,0)})$ for all $i$, and  $\height(I_{c,(n,d,t)})=\height(I_{c,n-(d-1)t,d,0)})$ by \cite[Proposition~3.7 (a)]{DHQ}.  On the other hand, we know that  $I_{c,n-(d-1)t,d,0)}=I_{\cb,n-(d-1)t,d}=I_{\cb,3,d}$ where $\cb=(c,\ldots , c) \in \ZZ^n$.  Hence, $\dim (I_{\cb,3,d})=2$. Since $a_1 \geq a_2 \geq a_3$,  it follwos that   $I_{\cb,n-(d-1)t,d} = x_1^{d'}J$ for a positive integer $d'$ and an ideal of Veronese type  $J$ with $\dim(J)=1$.  Set $\delta=d-d'$. So, we obtain the following corollary from Proposition~\ref{jjj}.
\begin{Corollary}
\label{nnn}
Let $I=I_{c,(n,d,t)}$ be a $c$-bounded $t$-spread Veronese ideal such that $n-(d-1)t=3$.
 If $\dim(I)=0$, then
\[
\beta_2(I)=d(d+2),  \quad  \beta_3(I)={d+1 \choose 2}.
\]
If $\dim(I)=1$, then
\[
\beta_2(I)=2\beta_1(I)-d-2,  \quad  \beta_3(I)=\beta_1(I)-d-1.
\]
If $\dim(I)=2$, then 
\[
\beta_2(I)=\beta_1(I)-\delta-2,  \quad  \beta_3(I)=\beta_1(I)-\delta-1.
\]
\end{Corollary}
 We also obtain the following corollary from Corollary~\ref{nnn} and the fact that $I_{n,d,t}=I_{d,(n,d,t)}$.
\begin{Corollary}
\label{hhh}
Let $I=I_{n,d,t}$ be a $t$-spread Veronese ideal of degree $d$ such that $n-(d-1)t=3$.
If $\dim(I)=0$, then
\[
\beta_2(I)=d(d+2),  \quad  \beta_3(I)={d+1 \choose 2}.
\]
If $\dim(I)=1$, then 
\[
\beta_2(I)=2\beta_1(I)-d-2,  \quad  \beta_3(I)=\beta_1(I)-d-1.
\]
If $\dim(I)=2$, then 
\[
\beta_2(I)=\beta_1(I)-\delta-2,  \quad  \beta_3(I)=\beta_1(I)-\delta-1.
\]
\end{Corollary}


\end{document}